\documentclass{article}
\usepackage{float}
\usepackage{algorithmicx}
\usepackage{graphicx}
\usepackage{amsmath}
\usepackage{amssymb}
\usepackage{tabularx}
\usepackage[ngerman,english]{babel}
\usepackage{authblk}
\usepackage{geometry}
\title{Using a one-dimensional finite-element approximation of Webster's horn equation to estimate individual ear canal acoustic transfer from input impedances}
\author[1]{Nick Wulbusch}
\author[2]{Reinhild Roden}
\author[1]{Alexey Chernov}
\author[2,3]{Matthias Blau}
\affil[1]{Institut für Mathematik,  Carl-von-Ossietzky Universität, Oldenburg, Germany}
\affil[2]{Institut für Audiologie und Hörtechnik, Jade-Hochschule, Oldenburg, Germany}
\affil[3]{Cluster of Excellence "Hearing4All"}
\affil[]{\textit{nick.wulbusch@uni-oldenburg.de}}
\date{}              
\begin{document}
	\maketitle
	\begin{abstract}
In many applications, knowledge of the sound pressure transfer to the eardrum is important. The transfer is highly influenced by the shape of the ear canal and its acoustic properties, such as the acoustic impedance at the eardrum. Invasive procedures to measure the sound pressure at the eardrum are usually elaborate or costly. In this work, we propose a numerical method to estimate the transfer impedance at the eardrum given only input impedance measurements at the ear canal entrance by using one-dimensional first-order finite elements and Nelder-Mead optimization algorithm. Estimations on the area function of the ear canal and the acoustic impedance at the eardrum are achieved. Results are validated through numerical simulations on ten different ear canal geometries and three different acoustic impedances at the eardrum using synthetically generated data from three-dimensional finite element simulations.   
	\end{abstract}

	\section{Introduction\label{sec:1}}
	
	Humans have individually shaped ear canals which will individually influence the sound transfer from the ear canal entrance to the eardrum. In many applications such as hearing aids or insert earphones, one is interested in the sound pressure at the eardrum, however, it is often not practical to measure it directly. Therefore, it is desirable to be able to individually predict the eardrum sound pressure from measurements at the entrance to the ear canal. One possibility is to measure the acoustic impedance, i.e., the ratio of sound pressure and volume velocity, at the entrance of the ear canal. The goal is then to estimate the transfer impedance, i.e., the ratio of the eardrum sound pressure to the volume velocity at the entrance, because one can then use input and transfer impedances together with a source model to individually predict the eardrum pressure, see e.g. \cite{Blau.2010}.\\ 
	A useful one-dimensional equation to study this behavior is Webster's horn equation \cite{Webster.1919}, describing the sound pressure in a horn, i.e., a tube with varying diameter. In order to solve Webster's equation for the sound pressure at different positions along the ear canal, an area function which describes the area along the ear canal is needed. If three-dimensional ear canal geometries are given, a curved center axis and an associated area function can be computed, using for instance the method proposed by \cite{Stinson.1989}. In most cases, such data is, however, not available. Alternatively, area functions can be estimated from acoustic impedances or reflectances measured at the entrance of the individual ear canal.\\
	In a number of publications, the phase of the reflectance was used to estimate the transfer impedance and area function \cite{Hudde.1999, Blau.2010, Rasetshwane.2011, SankowskyRothe.2011, SankowskyRothe.2015, Souza.2014}. As an example, in \cite{Hudde.1999} the optimization was done by an electro-acoustic model and a gradient method to approximate the radius function for a measured (or simulated) input reflectance. Given reference pressure values at the eardrum, the pressure was transformed back to the entrance for individual ear canals. A chain matrix, optimized using reflectance data then lead to accurate results for the radius function and the pressure at the drum for frequencies up to 16\,kHz.\\
	Another approach using a reflectance-based method was considered in \cite{Rasetshwane.2011} where the aim was to reproduce the ear canal cross-sectional area.  Specifically, a time-domain reflectance was calculated and the inverse solution to the time-domain horn equation was formulated in terms of forward and backward traveling pressure waves to obtain the area function by using a finite difference approximation. The solution to the inverse problem was analyzed for infinite acoustic horns in \cite{Rasetshwane.2012}. In \cite{Rasetshwane.2011} the area function of the inverse solution was applied to measurements made in an ear simulator. The estimated ear canal area function was close to the one of the simulator. Furthermore, the authors estimated the ear canal area function from reflectance measurements on real subjects, assuming the eardrum to be rigid, although the true ear canal area functions were not known. The resulting ear canal area functions looked plausible and similar to other ear canal area functions from the literature. The transfer impedance, however, was not investigated.\\
	In this paper, we present a different approach based on data of the input impedance in the frequency domain. Working in the frequency domain appears to be more convenient with respect to actual measurement data (i.e., acoustic impedances or reflectances) whose validity will in general be limited to a certain frequency range.  Our aim is to estimate the transfer impedance of the individual residual ear canal. To this end, we parameterize the ear canal area function and the acoustic impedance of the eardrum. We solve the horn equation in the frequency domain with linear finite elements and optimize the parameters using the Nelder-Mead method \cite{Nelder.1965} with respect to given data of the input impedance. By doing so, we also obtain estimations of the ear canal area function and the acoustic impedance at the eardrum.\\
	The general idea of the method is described in section \ref{sec:2}. The generation of data as input and for validation purposes is explained in section \ref{sec:3}. In  section \ref{sec:4} the procedure from section \ref{sec:2} will be extended by several adjustments to avoid unrealistic results during the parameter fitting. The final version of our method is validated in section \ref{sec:5} on a variety of ear canal geometries with different eardrum impedance models. Finally, in section \ref{sec:6} conclusions are given.\\
	A MATLAB implementation of the complete method can be found in the supplementary material \cite{Wulbusch}.\footnote{See supplementary
	material at https://github.com/nickwulbusch/ear-canal-parameter-fitting-1d for full MATLAB-Code of the final method.}
	
	\section{Method\label{sec:2}}
	The main goal of our method is to estimate the transfer impedance $Z_{\mathrm{tr}}$, based on data of the input impedance $Z_{\mathrm{in},\mathrm{data}}$, for individual ear canals. To do this, we use Webster's horn equation \cite{Webster.1919} in the frequency domain as surrogate model of a simplified representation of the pressure distribution in the ear canal,
	\begin{subequations}
		\begin{eqnarray}
			\frac{d}{dx}\left[S(x)\frac{d}{dx}p(x)\right]+k^{2}S(x)p(x) &=& 0  \text{ in } [0,\ell], \label{horneq}
			\\
			S(0)\frac{dp}{dx}(0) &=& q\mathrm{i}\omega \rho, \label{horneq0}
			\\
			S(\ell)\frac{dp}{dx}(\ell)+\frac{\mathrm{i}\omega\rho}{Z_{d}}p(\ell) &=& 0. \label{horneqL}
		\end{eqnarray}
	\end{subequations}
	The solution $p(x)$ to this problem denotes the acoustic pressure, $S(x)$ is the area function, $k=\frac{2\pi f}{\lambda}$ the wavenumber with frequency $f$ and wavelength $\lambda$, $q$ the volume velocity, $\omega=2\pi f$ the angular frequency, $\rho$ the density of air inside the ear canal and $Z_{d}$ the acoustic impedance of the eardrum. In this setting, the aim is to estimate the transfer impedance $Z_{\mathrm{tr}} = \frac{p(\ell)}{q}$ from given data of the input impedance $Z_{\mathrm{in},\mathrm{data}}$ which in turn is modeled by $Z_{\mathrm{in}} = \frac{p(0)}{q}$ in the one-dimensional surrogate \eqref{horneq}-\eqref{horneqL}. The area function $S(x)$ and acoustic impedance $Z_{d}$ at the eardrum are unknown. In the following, these quantities are parameterized and fitted to the given data of the input impedance by minimizing the cost function
	\begin{equation}
		\label{eq:J}
			J_{0}(Z_{\mathrm{in}}) = \sum\limits_{f\in\mathcal{F}} A\bigg\lvert \log_{10} \bigg\lvert \frac{Z_{\mathrm{in}}(f)}{Z_{\mathrm{in,data}}(f)}\bigg\rvert \bigg\rvert^{2} + B \arg\left(\frac{Z_{\mathrm{in}}(f)}{Z_{\mathrm{in,data}}(f)}\right)^{2},
	\end{equation}
	with weighting parameters $A$ and $B$. These are chosen as $A=10$ and $B=1$. The cost function is thus the weighted sum of the squared differences between model and data of the logarithmic amplitude and the phase in radians. $Z_{\mathrm{in,data}}(f)$ denotes the given data, i.e., the input impedance from the measurement at the entrance of the ear canal for a specific frequency $f$, whereas $Z_{\mathrm{in}}(f)$ denotes the input impedance from the solution of problem \eqref{horneq}-\eqref{horneqL} at the entrance. $\mathcal{F}$ denotes a set of frequencies. Different choices of $\mathcal{F}$ will be discussed in section \ref{subsec:highfrequency}. Unknown parameters are involved in the area function $S(x)$, the model for the acoustic impedance at the eardrum $Z_{d}$ and the length of the ear canal $\ell$. The area function is modeled as
	\begin{equation}
		\label{eq:area}
			S(x) := S(x,S_{0},\mathbf{c},\mathbf{s},\ell)= S_{0} + \sum\limits_{m=1}^{M} c_{m} \cos\left(\frac{m\pi x}{\ell}\right) + s_{m}\sin\left(\frac{m\pi x}{\ell}\right).
	\end{equation}
	The acoustic impedance $Z_{\text{ED}}$ at the eardrum is modeled using a two-resonator model, similar to the model described in \cite{A.Stirnemann.},
	\begin{equation}	\label{eq:Zd}
			Z_{\mathrm{ED}}= \bigg(\frac{1}{10^{L_{0,1}/20\,\text{dB}}(\mathrm{i}v_{1}Q_{1}+1)}+\frac{1}{10^{L_{0,2}/20\,\text{dB}}(\mathrm{i} v_{2}Q_{2}+1)}\bigg)^{-1}\,\text{Pa$\cdot$s/m$^{3}$},
	\end{equation}
	where
	\begin{equation}
		\label{eq:v}
		v_i = \frac{\omega}{2\pi f_{0,i}}-\frac{2\pi f_{0,i}}{\omega}
	\end{equation}
	and the parameters $Q_{i}$,  $f_{0,i}$ and $L_{0,i}$ are the quality factor, the resonance frequency and the impedance level at resonance (in dB re 1 Pa$\cdot$s/m$^{3}$), respectively. This model is a good compromise between the number of parameters and approximation as compared to more sophisticated models like the model by \cite{Hudde.1998c}. It is notably capable of modeling the cadaver measurements of \cite{Rosowski.1990} with acceptable accuracy.\\
	Following the suggestion by \cite{Hudde.1998a}, the innermost part of the ear canal is modeled as a lumped compliance (of a right circular cone of 4\,mm length and 2.5\,mm radius) acting in parallel with the eardrum impedance,
	\begin{equation}
		\label{eq:ZDCR}
		Z_{d} = \frac{Z_{\text{ED}}Z_{\mathrm{vol}}}{Z_{\mathrm{ED}}+Z_{\mathrm{vol}}},
	\end{equation}
	where
	\begin{equation}
		\label{eq:Zvol}
		Z_{\mathrm{vol}} = \frac{\rho c^{2}}{\mathrm{i}\omega V},
	\end{equation}
	$Z_{\text{ED}}$ corresponds to the acoustic impedance at the eardrum, $c$ is the speed of sound and $V$ is the volume of the cone. This effectively removes the innermost part from the ear canal geometry.\\
	The parameters used in the definitions of the area function, i.e., $S(x)$ in \eqref{eq:area} and the acoustic impedance of the eardrum $Z_{d}$ in \eqref{eq:ZDCR}, were fitted using the Nelder-Mead constrained optimization procedure \cite{Nelder.1965,fminsearchbnd}. To this end, equations \eqref{horneq}-\eqref{horneqL} were solved numerically by a self-implemented finite-element code in \texttt{MATLAB} using linear basis functions and the Simpson quadrature rule to compute the integrals. For this, equation \eqref{horneq} was multiplied by a test function $u(x)$ and integrated over the interval $[0,\ell]$. Changing the signs and integrating by parts lead to the so called weak formulation
	\begin{equation}
		\label{eq:weakFormulation}
			\int_{0}^{\ell} S(x) p'(x) u'(x) \,dx + \frac{i\omega\rho}{Z_{d}}p(\ell)u(\ell)- k^{2} \int_{0}^{\ell} S(x) p(x) u(x) \,dx = -qi\omega \rho u(0).
	\end{equation}
	To discretize this problem we divided the interval into smaller intervals $[x_{n-1},x_{n}]$ with $x_{n} = \frac{n\ell}{N}$, $n=1,\dots,N$ for some $N\in \mathbb{N}$. These subintervals are called elements. Next we considered a basis $\{\varphi_{n}\}_{n}$ of continuous and piecewise linear functions, where $\varphi_{n}$ are hat functions on the corresponding subinterval $[x_{n},x_{n+1}]$ and 0 elsewhere:
	\begin{equation}\label{eq:hatFunction}
		\varphi_{n}(x) = \begin{cases}
			\frac{x-x_{n-1}}{x_{n}-x_{n-1}}, & \text{if } x\in \left[\left.x_{n-1}{,} x_{n}\right)\right., n=1,\dots,N,\\
			\frac{x_{n+1}-x}{x_{n+1}-x_{n}}, & \text{if } x \in [x_{n},x_{n+1}]{,} n=0,\dots,{N-1},\\
			0, & \mathrm{else}.	
		\end{cases}
	\end{equation}
		\begin{figure}[t]
		\centering
		\includegraphics[width=0.55\textwidth]{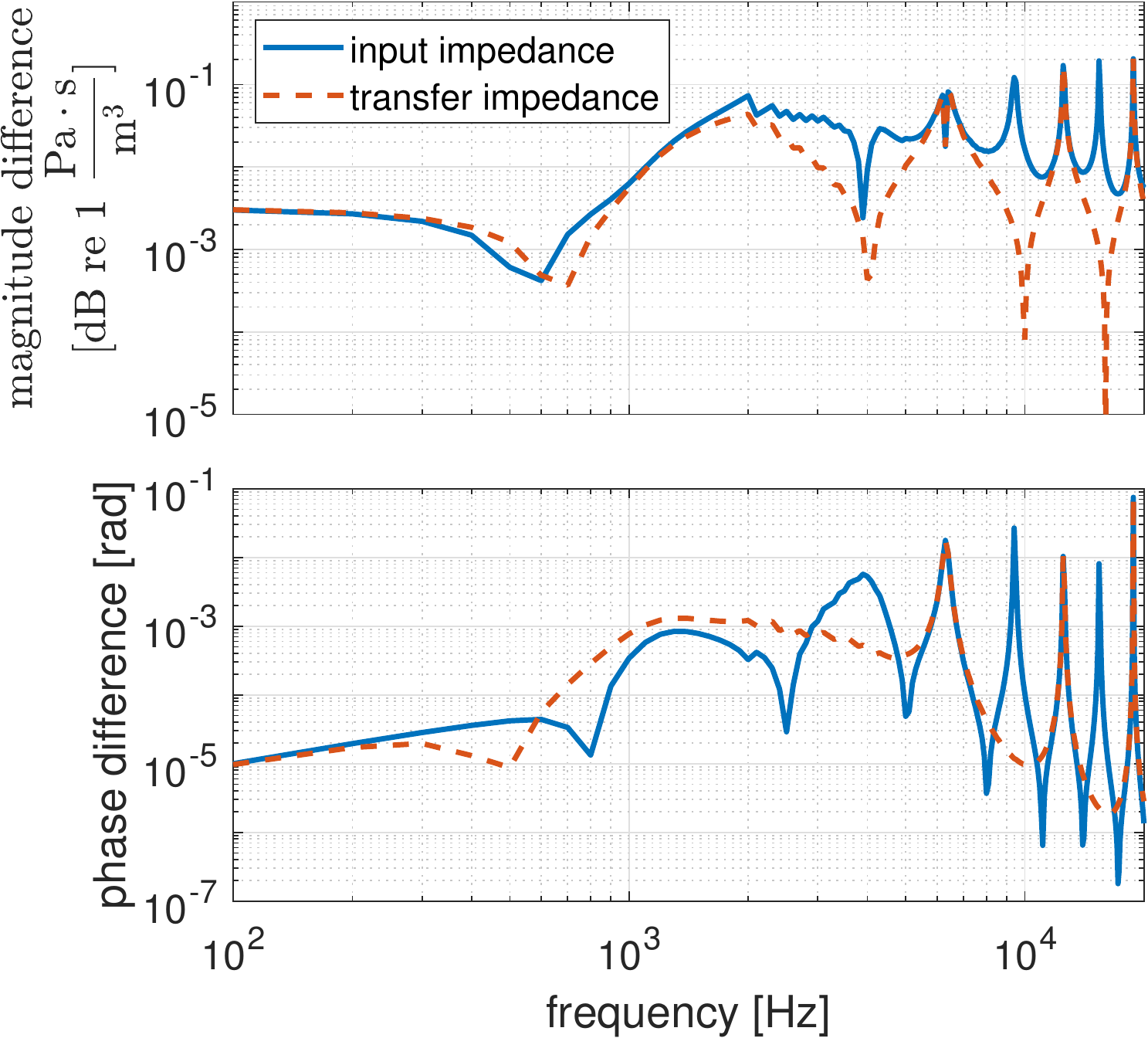}
		\caption{\label{FIG1}{Level difference and phase difference of $Z_{\mathrm{in}} = \frac{p(0)}{q}$ and $Z_{\mathrm{in,ref}} = \frac{p_{\mathrm{ref}}(0)}{q}$ and $Z_{\mathrm{tr}} = \frac{p(L)}{q}$ and $Z_{\mathrm{tr,ref}} = \frac{p_{\mathrm{ref}}(L)}{q}$. The discrete solution $p$ of \eqref{eq:weakFormulation} is computed with $N=4\max\left(1\,\mathrm{m}\cdot\ell/\lambda^{2},1\right)$ elements for a typical area function $S(x)$ and eardrum impedance $Z_{d}$. The reference solution $p_{\mathrm{ref}}$ is computed with $N_{\mathrm{ref}}=5000$ elements. The error is uniformly bounded within the frequency range 100\,Hz to 20\,kHz.}}
	\end{figure}
	The discrete solution $p_{N}$ to \eqref{eq:weakFormulation} is written in terms of this basis,
	\begin{equation}\label{eq:discreteSolution}
		p_{N}(x) = \sum\limits_{n=0}^{N} c_{n} \varphi_{n}(x),
	\end{equation}
	where the coefficients $c_{n}$ are computed by solving the system of linear equations
	\begin{equation}\label{eq:FEsystem}
		(A-M+R)c = F,
	\end{equation}
	where 
	\begin{subequations}\label{eq:FEmatrices}
		\begin{eqnarray}
			A_{ij} &=& \int_{0}^{\ell} S(x)\varphi_{j}'(x) {\varphi_{i}'(x)} \,dx,\\
			M_{ij} &=& k^{2} \int_{0}^{\ell} S(x) \varphi_{j}(x) {\varphi_{i}(x)} \,dx,\\
			R_{NN} &=& \frac{\mathrm{i}\omega\rho}{Z_{d}}, \quad R_{ij} = 0, \text{ for } i\ne N \ne j,\\
			F_{1} &=&  -q i\omega \rho, \quad F_{i} = 0, \text{ for } i\ne 0.
		\end{eqnarray}
	\end{subequations}
	In the experiments, the number of elements was chosen frequency-dependent as $N=4\max\left(1\,\mathrm{m}\cdot\ell/\lambda^{2},1\right)$, i.e., the discretization is different for each frequency. In this case the typical "rule of thumb" of 10 elements per wave length is guaranteed. Additionally for the case of the horn equation, the discretization error with respect to a the magnitude and phase seems to be uniformly bounded, see figure \ref{FIG1}. The discretization error was computed with respect to a reference solution with $N_{\mathrm{ref}}=5000$ elements.\\
	\begin{figure}[t]
	\centering
	\includegraphics[width=1\textwidth]{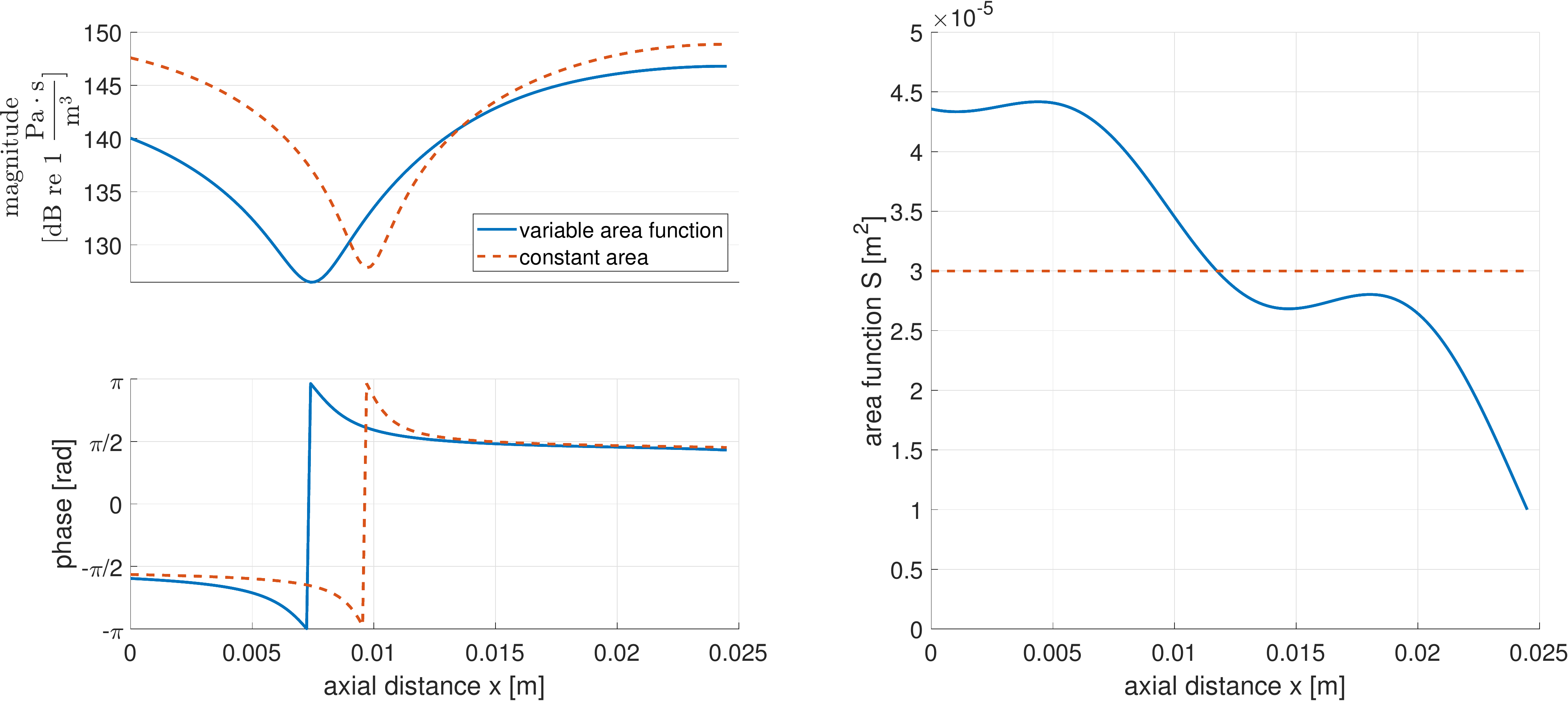}
	\caption{\label{FIG2}{Solution for the Horn equation with $f=6000\,$Hz with variable area function and constant area function using linear finite elements (top) and associated area functions (bottom).}}
\end{figure}
	Figure \ref{FIG2} illustrates an example solution for the weak formulation \eqref{eq:weakFormulation} of the horn equation \eqref{horneq}-\eqref{horneqL} for a frequency of $f=6\,$kHz, $Z_{d} = (6.08+7.34\mathrm{i})\times 10^{7}$\,$\text{Pa$\cdot$s/m$^{3}$}$ and an area function chosen as illustrated in figure \ref{FIG2}, i.e., a typical area function for an ear canal. In comparison, the solution to the pure Helmholtz problem, i.e., using a constant area function, is visualized as well to see the impact of the area function.

	\section{Data\label{sec:3}}
\begin{figure}[t]
			\centering
		\includegraphics[width=0.55\textwidth]{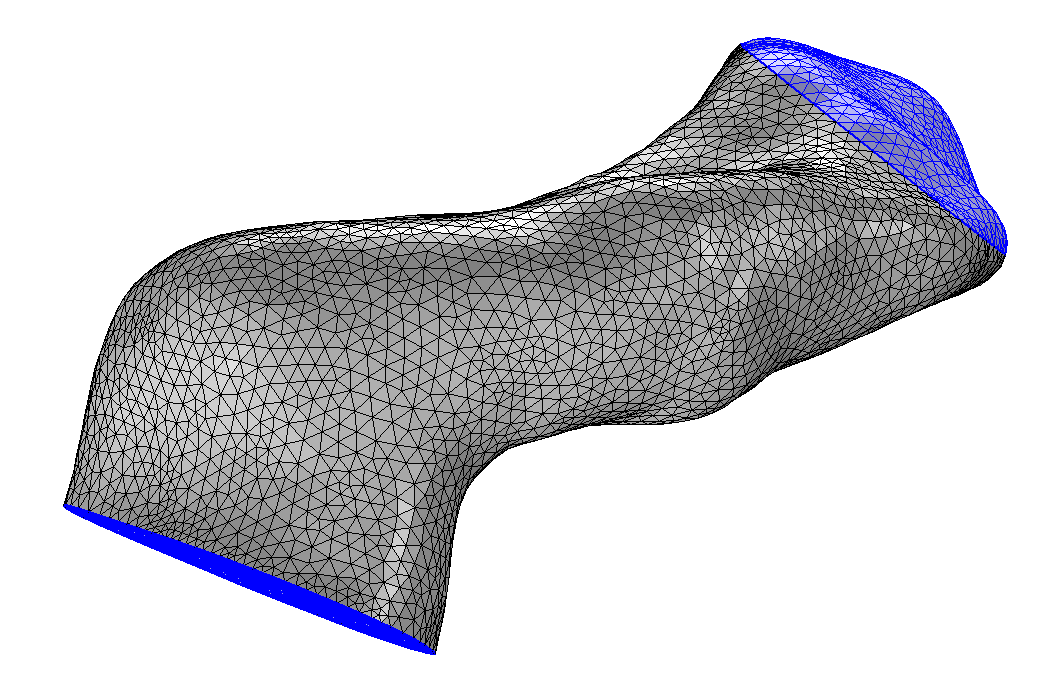}
		\caption{\label{FIG3}{Mesh of the ear canal of subject 5 cut at the first bend. The ear canal was rotated such that both the entrance and the eardrum (marked faces) are visible.}}	
	\end{figure}

	The synthetically generated input impedance and validation data were created by 3D finite-element simulations on geometries from the IHA database \cite{inproceedings,Roden.2021}. For the simulation, the STL geometry was cut at the first bend of the ear canal. The first bend was determined by computing the center axis of the ear canal using the VMTK toolbox \cite{vmtk}. Then, the coordinates with highest curvature of this center axis in the region of the first bend, which can be visually estimated, was evaluated using the \texttt{MATLAB} function of \cite{matlabCurve}. The implementation by \cite{frenet} was used to shape the entrance surface with respect to the coordinates of the first bend and the corresponding tangential vectors. This modified geometry was used in \cite{comsol}, where the geometries were discretized in a tetrahedral mesh with approximately 70{,}000-100{,}000 degrees of freedom, which corresponds to a maximum edge length of 1\,mm of an element. The meshes at the eardrum were refined further with maximum edge length of 0.2\,mm. As example, the meshed ear canal from subject 5 is illustrated in figure \ref{FIG3}. 

\begin{figure}[t]
			\centering
		\includegraphics[width=0.5\textwidth]{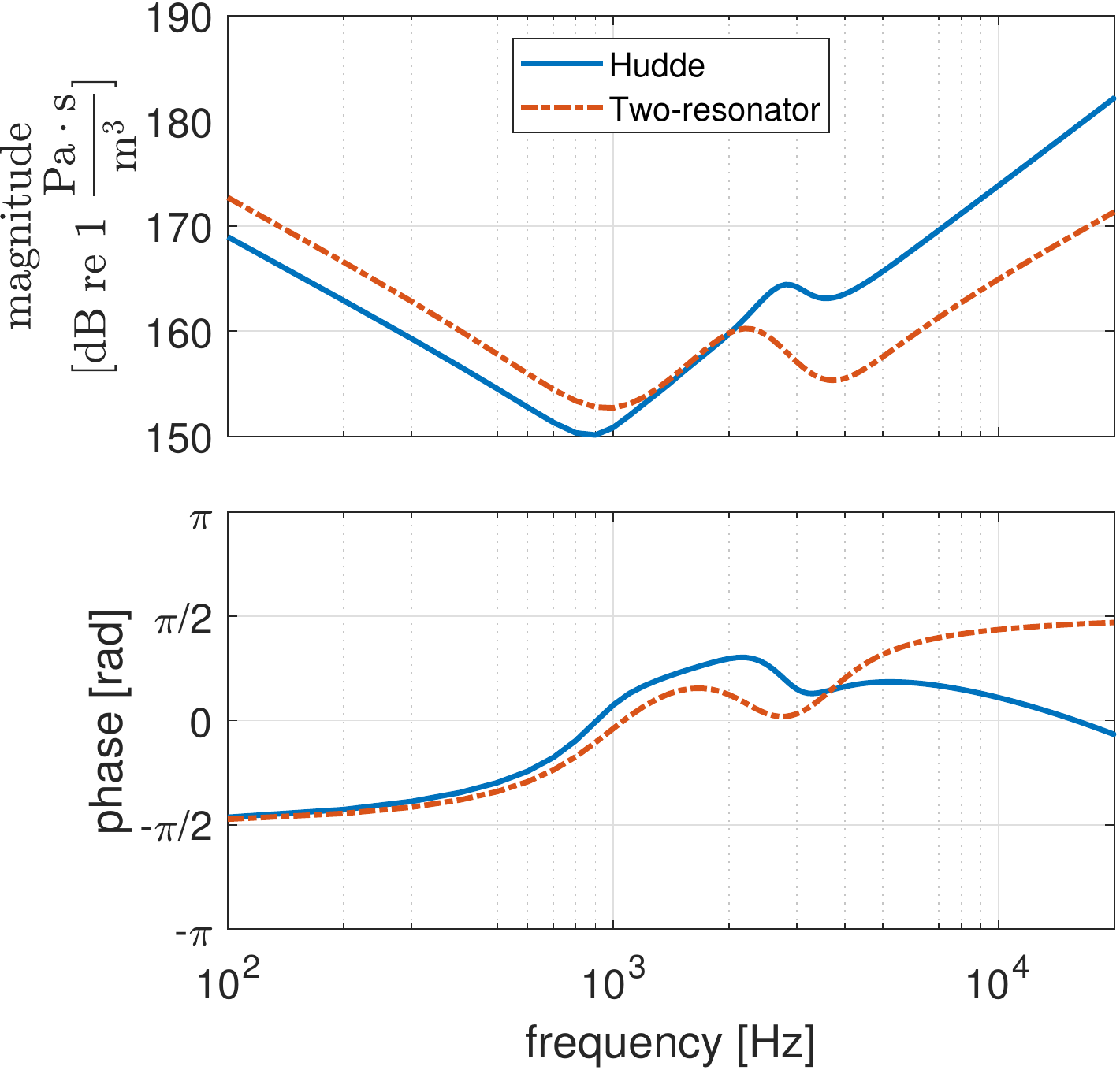}
		\caption{\label{FIG4}{Eardrum impedances types (i) and (ii) for $\widetilde{Z}_{ED}$ that were used to compute the reference data. Type (iii) is not shown due to the large magnitude that is constant at approximately 458\,dB with zero phase.}}
	\end{figure}

	The entrance surface was then assumed to vibrate piston-like to provide the excitation of the residual ear canal. Using this set-up, acoustic input impedances for three different types of impedances at the eardrum were computed at the entrance, by solving the Helmholtz equation,
	\begin{subequations}
			\begin{alignat}{4}
				\Delta p_{\mathrm{3D}}(\mathbf{x})+k^{2}p_{\mathrm{3D}}(\mathbf{x}) &= 0   &&\quad\text{ in } D, \label{eq:Helmholtz3D}
				\\
				\frac{\partial p_{\mathrm{3D}}}{\partial n} &= -v_{n}\mathrm{i}\omega \rho &&\quad\text{ on } \Gamma_{\text{entrance}},\label{eq:Helmholtz3Dentrance}
				\\
				\frac{\partial p_{\mathrm{3D}}}{\partial n}+\frac{\mathrm{i}\omega\rho}{\widetilde{Z}_{\mathrm{ED}}\lvert \Gamma_{\text{drum}}\rvert}p_{\mathrm{3D}} &= 0 &&\quad\text{ on } \Gamma_{\text{drum}}, \label{eq:Helmholtz3Ddrum}\\
				\frac{\partial p_{\mathrm{3D}}}{\partial n} &= 0 &&\quad\text{ on } \Gamma_{\text{wall}}, \label{eq:Helmholtz3Dwall}
		\end{alignat}
	\end{subequations}
	with $D$ being the interior of the geometry, $\Gamma_{\text{entrance}} \subset \partial D$ is the surface at the entrance, $\Gamma_{\text{drum}} \subset \partial D$ is the surface that corresponds to the eardrum and $\Gamma_{\text{wall}}$ is the wall of the ear canal. $\frac{\partial}{\partial n}$ denotes the outer normal derivative. Note that comparing equation \eqref{horneq0} and \eqref{eq:Helmholtz3Dentrance} we have different signs, since in \eqref{horneq0} the derivative in the one-dimensional model is directed into the interval. $v_{n}$ is the particle velocity. The input impedance used as input data during the parameter fitting was averaged over the entrance surface $\Gamma_{\text{entrance}}$
	\[
	Z_{\mathrm{in},\mathrm{data}} = \frac{\frac{1}{\lvert \Gamma_{\text{entrance}}\rvert}\int_{\Gamma_{\text{entrance}}} p_{\mathrm{3D}}(\mathbf{x}) d\Gamma_{\text{entrance}}}{\lvert \Gamma_{\text{entrance}}\rvert \text{ } v_{n}}.
	\] 
	The considered eardrum impedances $\widetilde{Z}_{\mathrm{ED}}$ are of the following types
	\begin{itemize}
		\item[(i)] impedance model by Hudde and Engel \cite{Hudde.1998c},
		\item[(ii)] a two-resonator model as described in \eqref{eq:Zd} with parameters given in table \ref{tab:table1},
		\item[(iii)] (nearly) rigid, i.e., $\widetilde{Z}_{\mathrm{ED}}\approx 8.4\times 10^{22} - 8.8 \times 10^{15}\mathrm{i}$\,$\text{Pa$\cdot$s/m$^{3}$}$,
	\end{itemize}
	They are illustrated in figure \ref{FIG4}. For the validation in section \ref{sec:5}, the transfer impedance $Z_{\mathrm{tr},\mathrm{data}}$ was computed with respect to the umbo, i.e.,
	\[
	Z_{\mathrm{tr},\mathrm{data}} = \frac{p_{\mathrm{3D}}(\mathbf{x}_\mathrm{umbo})}{\lvert \Gamma_{\text{entrance}}\rvert v_n}.
	\]
	The entrance impedance $Z_{\mathrm{in},\mathrm{data}}$ and transfer impedance $Z_{\mathrm{tr},\mathrm{data}}$ were computed for 200 distinct frequencies linearly distributed in the range between 100\,Hz and 20\,kHz. A total of 30 data sets, ten subjects with three different impedance models each, were available for testing and validation.
	\begin{table}[H]
	\centering
			\begin{tabular}{c|cc}
				& resonator 1 & resonator 2\\
				\hline
				$f_{0}$& 1000\,Hz & 3500\,Hz\\
				Q& 1.1 & 1.5\\
				$L_{0}$& 153\,dB re 1 $\mathrm{Pa}\cdot\mathrm{s}/\mathrm{m}^{3}$ & 157\,dB re 1 $\mathrm{Pa}\cdot\mathrm{s}/\mathrm{m}^{3}$
			\end{tabular}
		\caption{\label{tab:table1}Parameters used for the generation of data using the two-resonator model.}
	\end{table}
	\section{Original parameter fitting and refinement of the method\label{sec:4}}
\begin{figure}[t]
		\includegraphics[width=\textwidth]{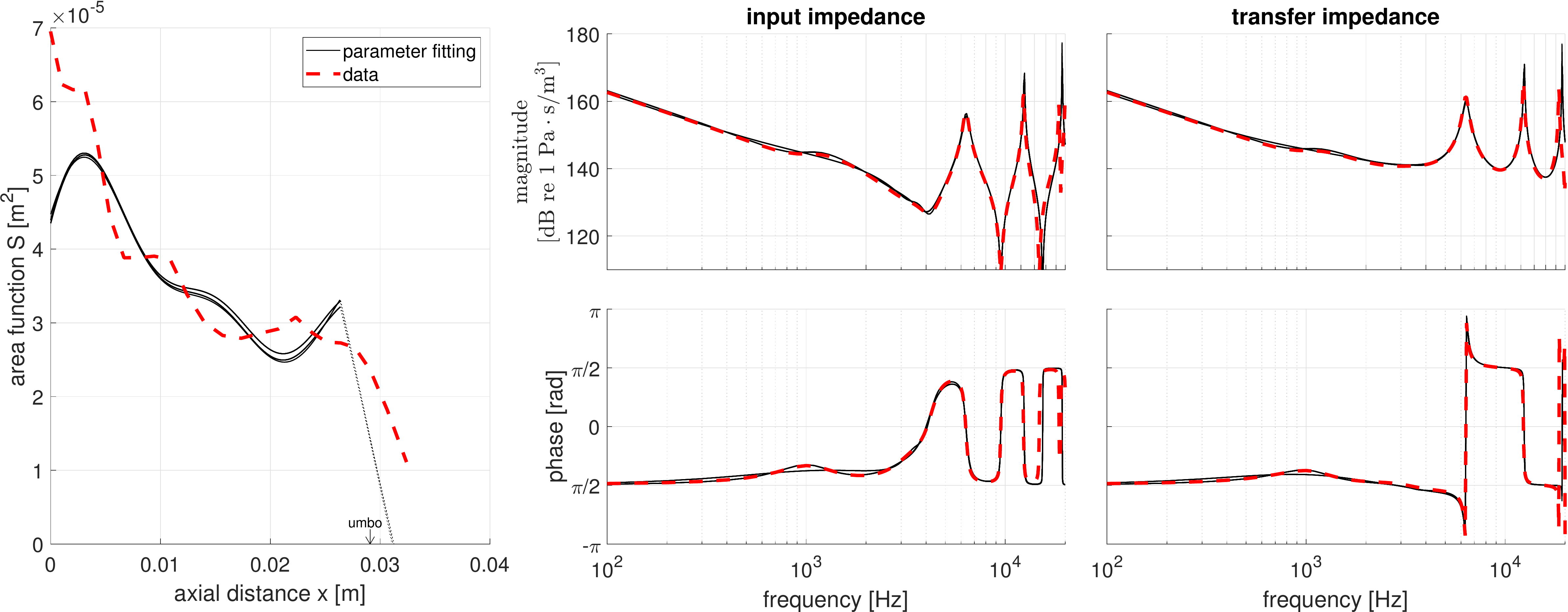}
		\caption{\label{FIG5}{Results for the original parameter fitting method for subject 5 and impedance model (ii). Solutions to each of the twelve different initial parameter sets are depicted by a black line respectively. Left: Area function. The rear part represents the cone as part of the impedance model, see discussion in section \ref{sec:2}. Center: Input impedance $Z_{\mathrm{in}}$. Right: Transfer impedance $Z_{\mathrm{tr}}$.}}
	\end{figure}

	\begin{table}[hbt!]
		\centering
		\small{
		\begin{tabular}{c|cccc}
			&$S_0$[$\text{m}^{2}$] &$\ell$[m] &$c_m$&$s_m$ \\
			\hline
			upper bound\rule{0pt}{2.6ex}& $2\times10^{-4}$ & $15\times 10^{-3}$ & $2^{-m+2}\times10^{-5}$ &$2^{-m+2}\times10^{-5}$\\
			lower bound& $10^{-5}$ & $45\times 10^{-3}$ & $-2^{-m+2}\times10^{-5}$ &$-2^{-m+2}\times10^{-5}$\\
			basic & $6\times 10^{-5}$ & $30\times 10^{-3}$ & $(2 \times 10^{-6},0,\dots,0)$ & $(0,0,\dots,0)$
		\end{tabular}\\[2ex]
		\begin{tabular}{c|ccccccccccccc}
			&$L_{0,1}$[dB] &$L_{0,2}$[dB] & $Q_{1}$ & $Q_{2}$ & $f_{0,1}$[Hz] & $f_{0,2}$[Hz] & $V$[$\text{m}^{3}$]\\[1ex]
			\hline
			upper bound\rule{0pt}{2.6ex}& 200  & $L_{0,1}+$40 & 10 & 10 & 2500 & 6000& $5.23\times 10^{-8}$\\
			lower bound& 50 & $L_{0,1}+$0 & 0.3 & 0.3 & 500 & 2500 & $1.3\times 10^{-8}$ \\
			basic & 161 & $L_{0,1}+$20 & 1.2 & 1.2 & 900 & 4000 & $2.62\times 10^{-8}$
	\end{tabular}}
	\caption{\label{tab:table2} Basic estimate with lower and upper bounds of the parameters of the ear canal area function and the two-resonator impedance model used in the one-dimensional surrogate model. The decibel values for impedance levels $L_{0,1}$ and $L_{0,2}$ refers to dB re 1 $\mathrm{Pa}\cdot\mathrm{s}/\mathrm{m}^{3}$.}
	\end{table}
	We began with the basic algorithm and $M=4$, i.e., four sine and cosine summands in the definition of the area function $S(x)$ in \eqref{eq:area}. The  parameter fitting was done with a subset of all available frequencies, more specifically with 25 logarithmically distributed frequencies in the range from 100\,Hz to 20\,kHz. The parameter constraints were set according to table \ref{tab:table2} to ensure that the area function $S(x)$ does not take excessive values, the oscillations of the area function are not unreasonably large and the impedance magnitude at the eardrum corresponds approximately to the expectations from commonly used models (i) and (ii) as in figure \ref{FIG4}. To enhance the accuracy, the Nelder-Mead method was used with three restarts, taking the result of the preceding optimization as initialization. Since the Nelder-Mead optimization procedure does not guarantee convergence to a global minimum, the  parameter fitting procedure was executed several times with different initial parameters. These were chosen as random perturbation of up to 25\,\% from the basic set of parameters, also depicted in table \ref{tab:table2}, i.e., for a given parameter $\alpha$ with basic value $\alpha_{0}$,
	\[
	\alpha_{\mathrm{init}} = \alpha_{0} \left( 1 + 0.25\, \mathcal{U}(-1,1)\right),
	\]
	where $\mathcal{U}(-1,1)$ is the uniform distribution over the interval $(-1,1)$. In total, twelve different initial parameter sets were considered. Figure \ref{FIG5} shows the results obtained with all twelve sets of initial parameters, considering data generated from the geometry of subject 5 and eardrum impedance model (ii). All twelve initial parameter sets converged to similar results, thus the lines for each model type are very close to each other. The area function $S(x)$ has a similar trend in comparison to the area function computed by the method of \cite{Stinson.1989} which is depicted as data curve in figure \ref{FIG5}. The input impedance is fitted well over the whole frequency range and even for the transfer impedance only slight differences between data and model can be observed.\\
	However, this original parameter fitting did not lead to comparably good results for all subjects and all different types of eardrum impedances. To enhance the robustness of the method, several aspects that lead to problems or unrealistic results for individual subjects and eardrum impedance models are considered in the following and possible solution strategies are discussed, implemented and illustrated.

	\subsection{Avoiding negative areas}
\begin{figure*}[t]
		\includegraphics[width=1\textwidth]{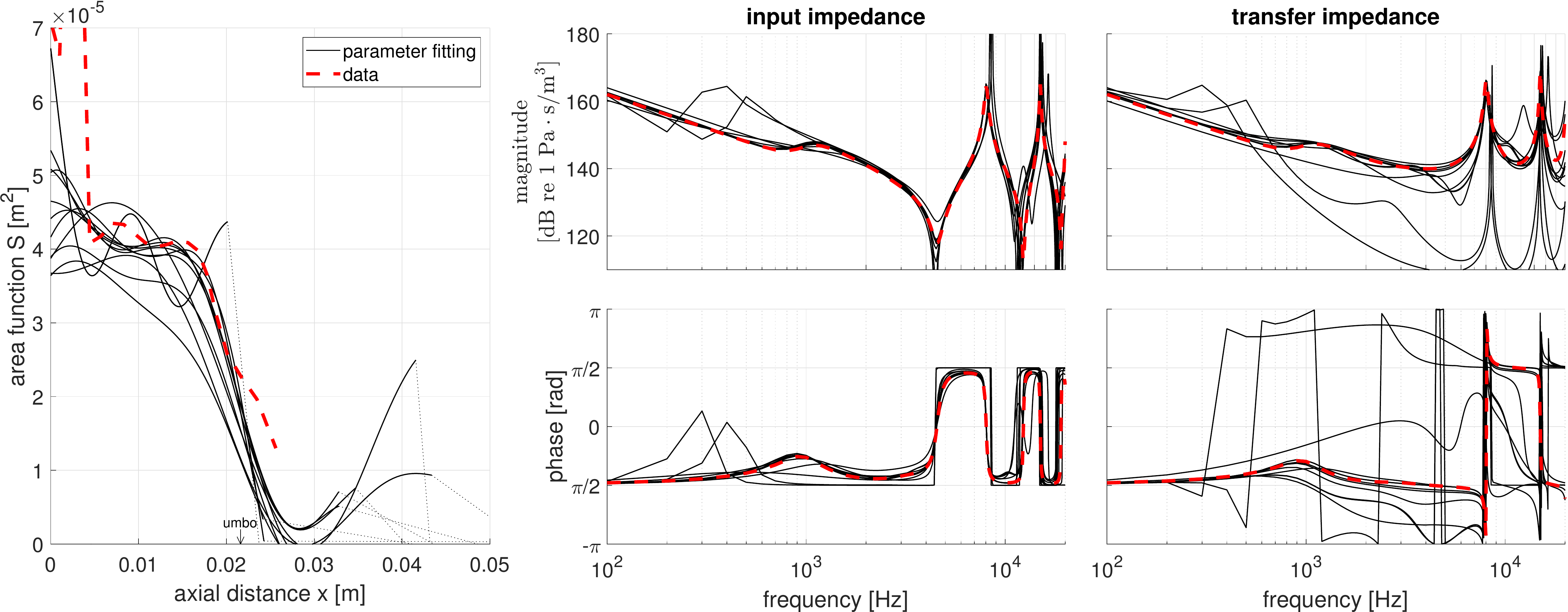}
		\caption{\label{FIG6}{Parameter fitting may lead to negative values of the area function. Results for subject 3 and impedance model (i). Solutions to each of the twelve different initial parameter sets are depicted by a black line respectively. Left: Area function. The rear part represents the cone as part of the impedance model, see discussion in section \ref{sec:2}. Center: Input impedance $Z_{\mathrm{in}}$. Right: Transfer impedance $Z_{\mathrm{tr}}$.}}
	\end{figure*}

	For some initial parameter sets, the original parameter fitting procedure returns a parameter set where the corresponding area function $S(x)$ takes negative values for some part in the interval $[0,\ell]$. Even though the reference data, i.e., the input impedance, was matched reasonably well, see figure \ref{FIG6}, some larger differences can be observed for the transfer impedance. We also note that the initial parameter set largely influences the result, thus the method is not very robust with regard to the initialization. Two possible strategies to prevent negative area function are
	\begin{enumerate}
		\item[(1)] further restriction on the parameters $S_{0},c_{m}$ and $s_{m}$ such that the area function $S$ will always be positive,
		\item[(2)] adding a penalty term to the cost function that penalizes negative or small values of $S$.
	\end{enumerate}
	Here the first strategy is not very practical. It leads to strong restrictions on $S_{0}$ in the way that the lower bound needs to be relatively large or the restrictions on $c_{m}$ and $s_{m}$ are such that admissible geometries are quite restricted.\\
	The second strategy still keeps a lot of freedom in the geometry-related parameters and will be used in the following. The penalty term is chosen as
	\begin{equation}
		J_{1,H^{(1)}}(S_{0},\mathbf{c},\mathbf{s},\ell) = 10^{4} \frac{\lVert \max\left(H^{(1)}-S(x),0\right)\rVert_{\infty}}{H^{(1)}},
	\end{equation}
	where $\lVert \cdot \rVert_{\infty}$ is the discrete supremum norm. This penalty term penalizes the area function becoming smaller than $H^{(1)}$. We choose $H^{(1)}=10^{-5}$\,m$^{2}$, which is about one half of the base area of the standard cone with a radius of 2.5\,mm.\\
	The updated procedure now usually generates functions that are positive on $[0,\ell]$. If the area still becomes negative then the model-data misfit is probably just too large and a meaningful parameter fitting is not possible with the underlying one-dimensional surrogate model \eqref{horneq}-\eqref{horneqL} or the initial values were chosen poorly.\\
		
	\begin{figure}[hbt!]
		\includegraphics[width=\textwidth]{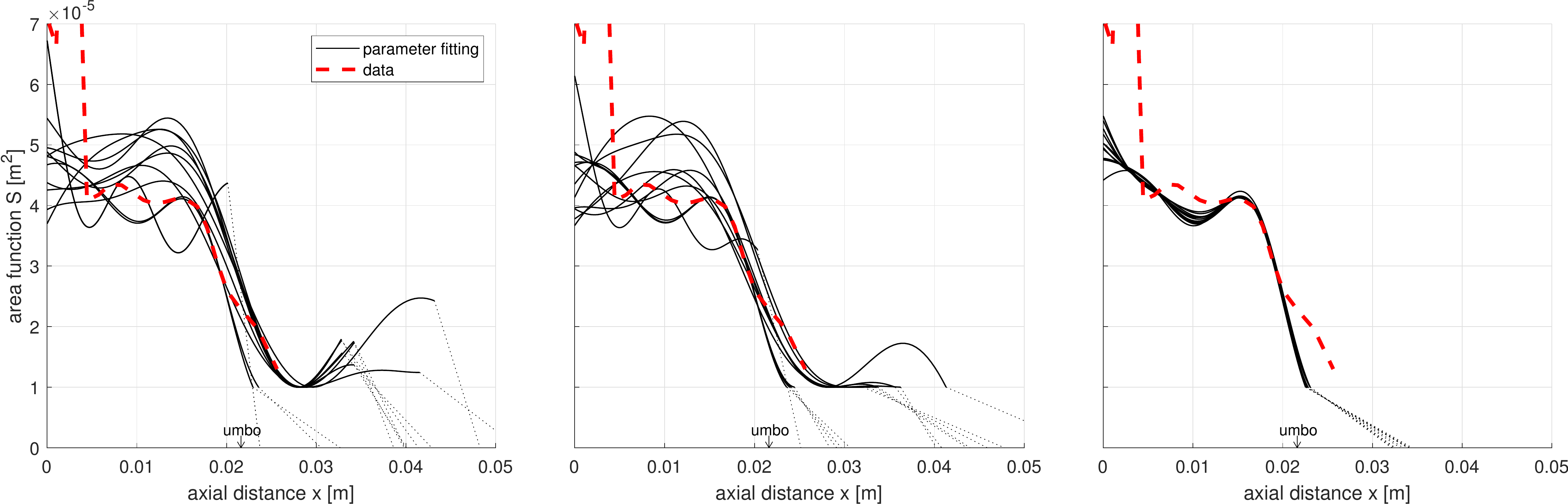}
		\caption{\label{FIG7}{Left: Results incorporating penalty term $J_{1}$ for subject 3 and impedance model (i) for twelve different initial parameter sets.  Parameter fitting may lead to unrealistic ear canal geometries due to large increases in the area function near the eardrum. Center: Results incorporating penalty terms $J_{1}$ and $J_{2}$ for subject 3 and impedance model (i) for twelve different initial parameter sets.  Parameter fitting may lead to largely overestimating the length of the ear canal. Right: Results incorporating penalty terms $J_{1}$, $J_{2}$ and prior estimating of the ear canal length $\ell$ for subject 3 and impedance model (i) for twelve different initial parameter sets.  Parameter fitting leads to a realistic are function.}}
	\end{figure}
	Occasionally the problem arises that the area function steeply increases after reaching the minimum. One case for this behavior is illustrated in figure \ref{FIG7}a. This is anatomically at least unusual. Therefore, another penalty term is applied that penalizes solutions in which the end of the ear canal is larger than the overall minimum or close to that minimum. This additional penalty term is given as
	\begin{equation}
			J_{2,H^{(2)}}(S_{0},\mathbf{c},\mathbf{s},\ell) = 10^4\left\lVert \max\limits \left( \frac{S(\ell)-S(x)-H^{(2)}}{\lvert S(x)\rvert},0\right) \right\rVert_{\infty}.
	\end{equation}
	The penalty is chosen as $H^{(2)}=0$ which enforces the minimum of the area function to be at the end of the ear canal. The cost function now becomes
	\begin{equation}
		\label{eq:J2Pen}
		J = J_{0} + J_{1,H^{(1)}} + J_{2,H^{(2)}}.
	\end{equation}
	Using the updated cost function $J$ in the parameter fitting leads to more realistic ear canal area functions, see figure \ref{FIG7}b.
	\subsection{Accomodating long and short ear canals}
	The length of the ear canal varies for each individual, thus the bounds of $\ell$ are chosen as $15\,\text{mm} \le \ell \le 45\,\text{mm}$ to be able to model short and long ear canals. Arbitrary initial values for $\ell$ may, however, lead the Nelder-Mead algorithm converge to a local minimum with an inadequate result for the value of $\ell$, as can be seen in figure \ref{FIG7}b. To avoid this, we estimate $\ell$ in a first step by looking at the first maximum of the input impedance magnitude of the data, see for example figure \ref{FIG6}b at around 8\,kHz. Assuming a rigidly terminated cylindrical waveguide, the initial length is estimated as
	\begin{equation}
		\label{eq:Lest}
		\ell_{\text{est}} = \frac{c}{2 f_{\text{max}}},
	\end{equation}
	where $f_{\text{max}}$ corresponds to the frequency where the first maximum of the input impedance magnitude, i.e., the data, is achieved. It turned out that in general the first maximum frequency gave a better result than the first minimum, which is why the length was estimated using equation \eqref{eq:Lest}. The ear canal is not a rigidly terminated cylindrical waveguide, but still this is a reasonable first approximation. Using this estimated $\ell_{\text{est}}$ as initial parameter leads to more appropriate results in the numerical simulations. To ensure that the overall length does not vary much from this estimated initial length during the  parameter fitting, the bounds for $\ell$ were redefined to $\ell_{\text{est}}-3\text{\,mm}\le \ell \le \ell_{\text{est}}+1\text{\,mm}$. The upper bound is stricter since the impedance model already assumes that the ear canal is longer due to the cone at the end representing the innermost part of the ear canal. These bounds lead to reasonable results for most initial values, see figure \ref{FIG7}c. Further results using the refined method are discussed in the validation in section \ref{sec:5}.\\
	In the following subsection we further discuss the influence of the chosen frequency set for the parameter fitting and the influence of the length of the expansion of the area function.

	\subsection{Influence of the bandwidth and frequency spacing used in the  parameter fitting}
	\label{subsec:highfrequency}
\begin{figure}[hbt!]
			\centering
		\includegraphics[width=1\textwidth]{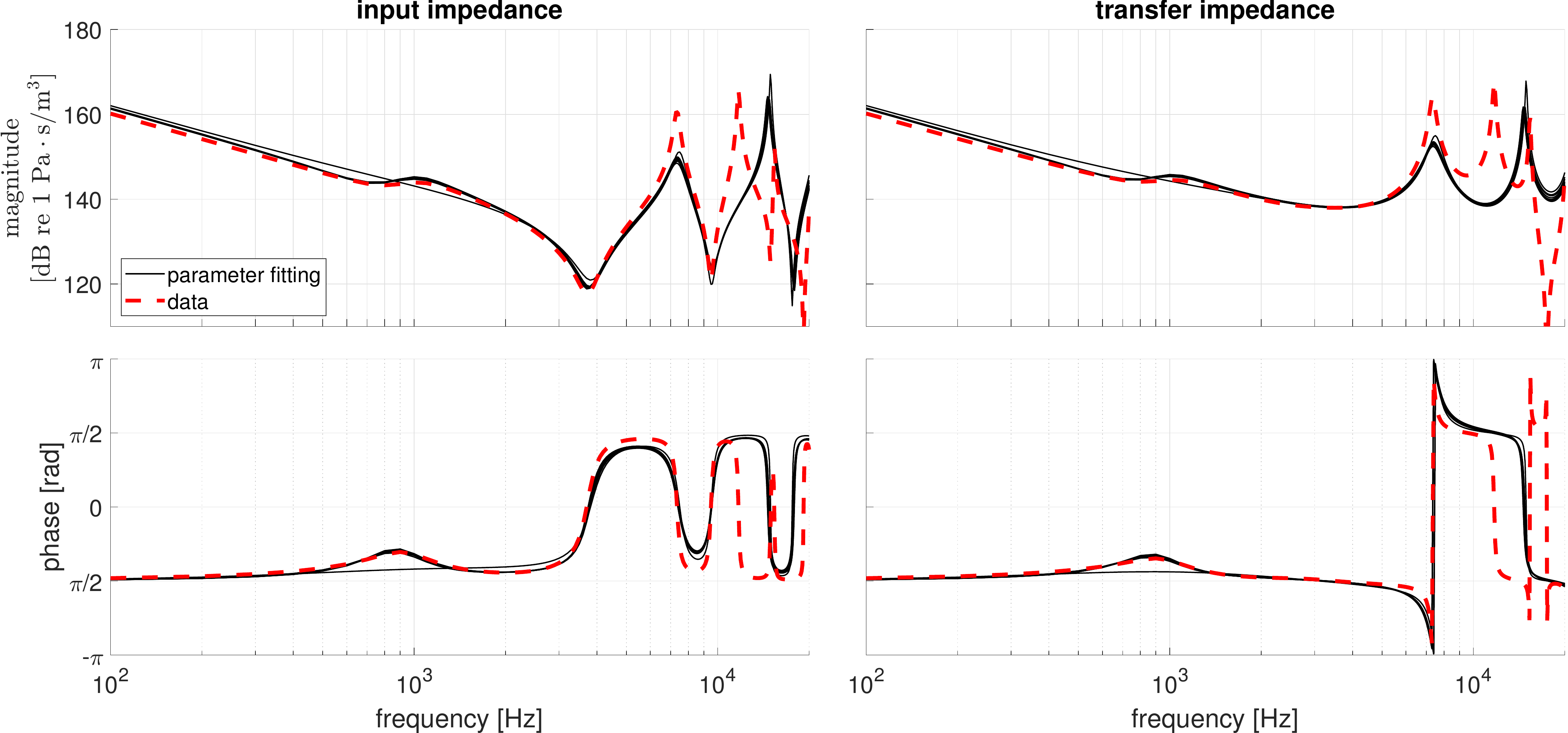}
		\caption{\label{FIG8}{Ear canal of subject 9 and impedance model (i) as example for a case of high frequency behavior that cannot be reproduced with the one-dimensional model. The last minimum in the entrance and transfer impedance data gets very low. Including high frequencies in the parameter fitting very likely influences precision for lower frequencies.}}
	\end{figure}

\begin{figure}[hbt!]
		\centering
		\includegraphics[width=1\textwidth]{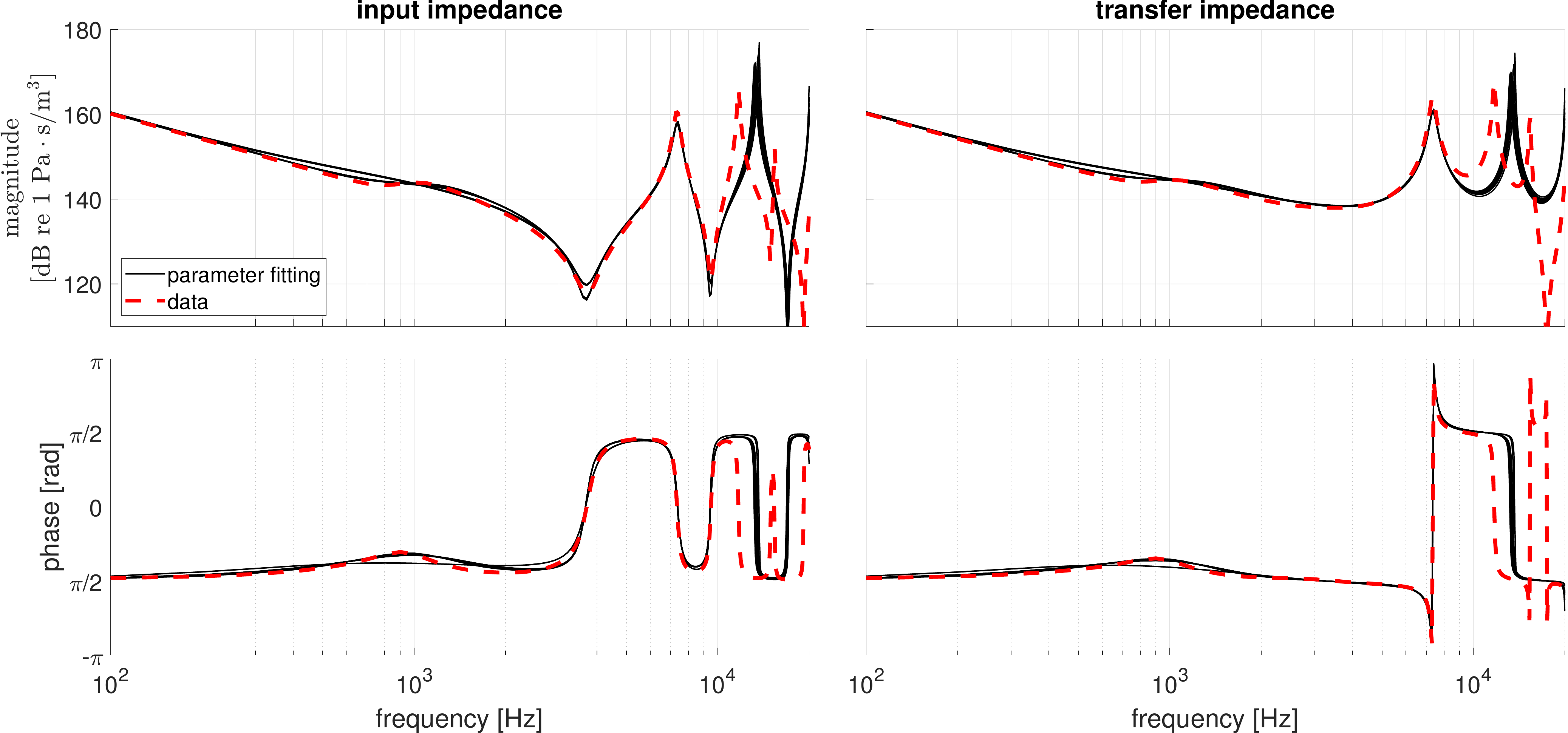}
		\caption{\label{FIG9}{Ear canal of subject 9 and impedance model (i) when only using data for frequencies up to 10\,kHz. Improvements in the lower frequency range, especially up to 8\,kHz can be observed in comparison to figure \ref{FIG8}.}}
	\end{figure}

	While the improvements of the prior subsections lead to better results in the  parameter fitting, in some cases problems arise at very high frequencies, because some features of the 3D geometry cannot be captured in a one-dimensional model. This means that not for all ear canal shapes good estimations on the transfer impedance for frequencies up to 20\,kHz can be expected. An example is illustrated in figure \ref{FIG8}. Here, the minimum of the magnitude of the input impedance at approximately 17\,kHz is much lower than usual. For this individual, good approximations for high frequencies cannot be achieved using the horn equation. It might therefore be useful to restrict the frequency range during the parameter fitting to get better results regarding the lower frequency part. Restricting the frequency range to frequencies smaller than 10\,kHz lead to improvements in the input impedance $Z_{\text{in}}$ and transfer impedance $Z_{\text{tr}}$ from the parameter fitting, see figure \ref{FIG9}.\\

\begin{figure}[hbt!]
		\centering
		\includegraphics[width=0.55\textwidth]{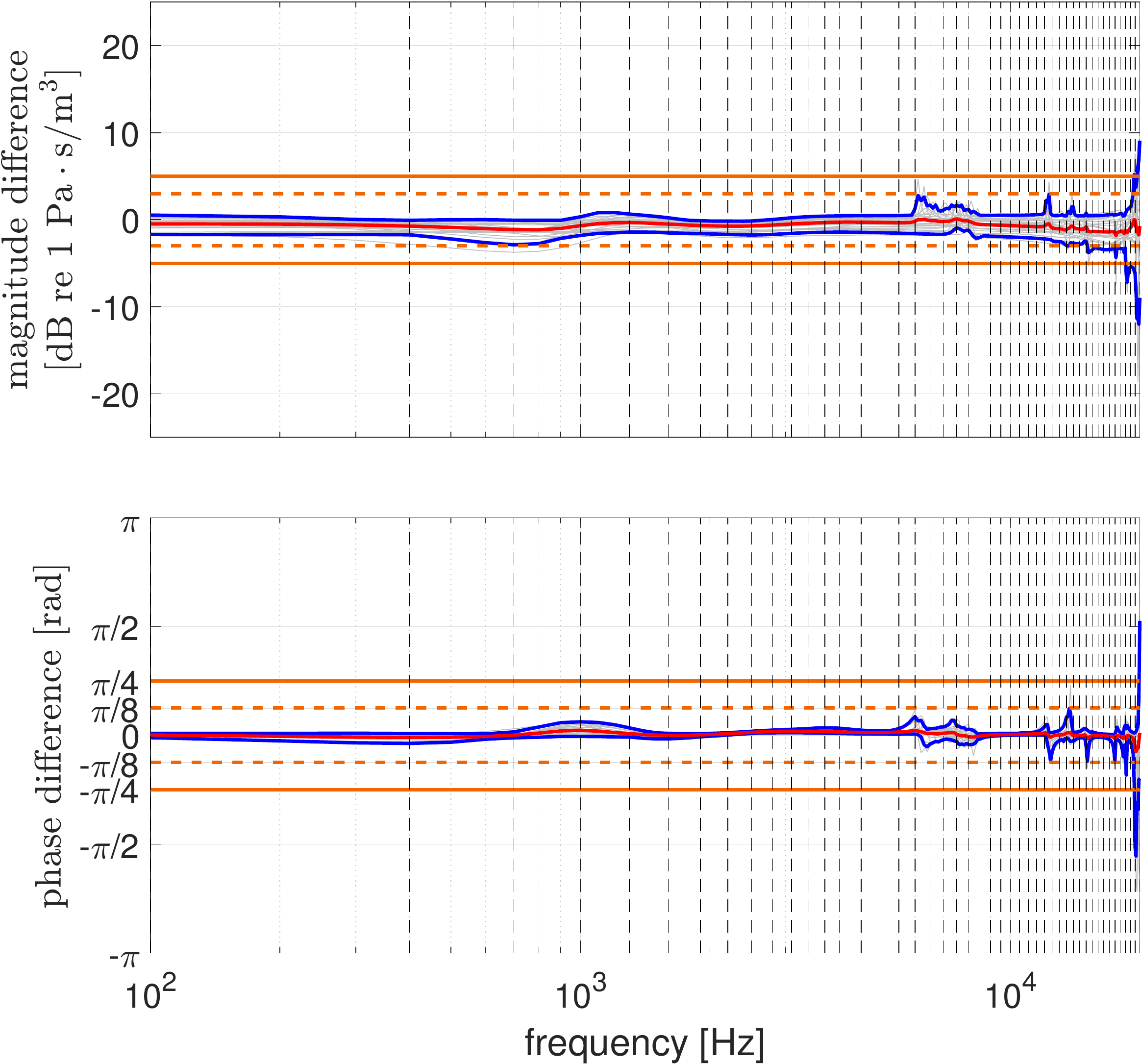}
		\caption{\label{FIG10}{Validation for the  parameter fitting considering differences of the transfer impedance for frequencies for linearly distributed frequencies up to 20\,kHz considering all data sets. Straight solid line: 5\,dB or 45° difference. Straight dashed line: 3\,dB or 22.5° difference. Blue lines: 5 and 95 quantile.  Red line in between: Mean. Dashed gray lines: Frequencies used during the parameter fitting.}}
	\end{figure}

	In the following, different choices for frequency subsets are considered. The considered frequency sets were either uniformly or logarithmically distributed (rounded up to the nearest 100\,Hz).  Figure \ref{FIG10} illustrates the validation for uniformly distributed frequencies where the magnitude and phase difference of the transfer impedance was computed with respect to the reference data. For some frequency sets, finer resolutions in the low frequency range (up to 3 kHz) or the higher frequency range (larger than 3 kHz) were used. Different upper limits were also considered, where the frequency sets were restricted either to frequencies up to 20\,kHz or 10\,kHz for reasons discussed above. Additionally, for some frequency sets $\mathcal{F}$, the frequencies that correspond to the maxima and minima of the input impedance magnitude were included, because these are the frequencies where usually the largest error with respect to the data is made and these frequencies are characteristic for the specific ear canal shape. We could observe that the frequency set using logarithmic distribution while additionally including the frequencies corresponding to the maxima and minima of the magnitude of the input impedance results in comparably low function values for the specific frequency range with respect to the validation function
	\begin{equation}
		\label{eq:Jval}
			J_{\mathrm{val}}(Z_{\mathrm{tr}}) = \sum\limits_{\substack{f=100\,\mathrm{Hz},200\,\mathrm{Hz},\\\dots,20\,\mathrm{kHz}}} A\bigg\lvert \log_{10} \bigg\lvert \frac{Z_{\mathrm{tr}}(f)}{Z_{\mathrm{tr,data}}(f)}\bigg\rvert \bigg\rvert^{2}+ B \arg\left(\frac{Z_{\mathrm{tr}}(f)}{Z_{\mathrm{tr,data}}(f)}\right)^{2}.
	\end{equation}
	Although other configurations lead to comparable validation function values, plotting the actual difference of the magnitude in dB over the frequency showed that for low frequencies the versions using uniformly distributed frequencies result in slightly worse behavior in the low frequency range with magnitude differences of up to 3\,dB in the frequency range from 100\,Hz to 6\,kHz, see figure \ref{FIG10} as example, while the logarithmically distributed frequencies only lead to magnitude differences of up to 1.5\,dB in this frequency range. This case is discussed and illustrated in the validation in section \ref{sec:5}. Consequently, a logarithmic distributed frequency set with maxima and minima frequencies was the choice for the validation in section \ref{sec:5}.\\

	Since usually frequencies smaller than 10\,kHz are especially of interest, we restrict the frequencies to 10\,kHz to improve the results in this range, especially if the pressure distribution in the ear canal exhibits large three-dimensional effects due to the geometry of the ear canal.

	\subsection{Influence of the number of parameters in the area function expansion} \label{subsec:FourierCoeff}
	 \begin{figure}[hbt!]
	 	\centering
	 	\includegraphics[width=0.55\textwidth]{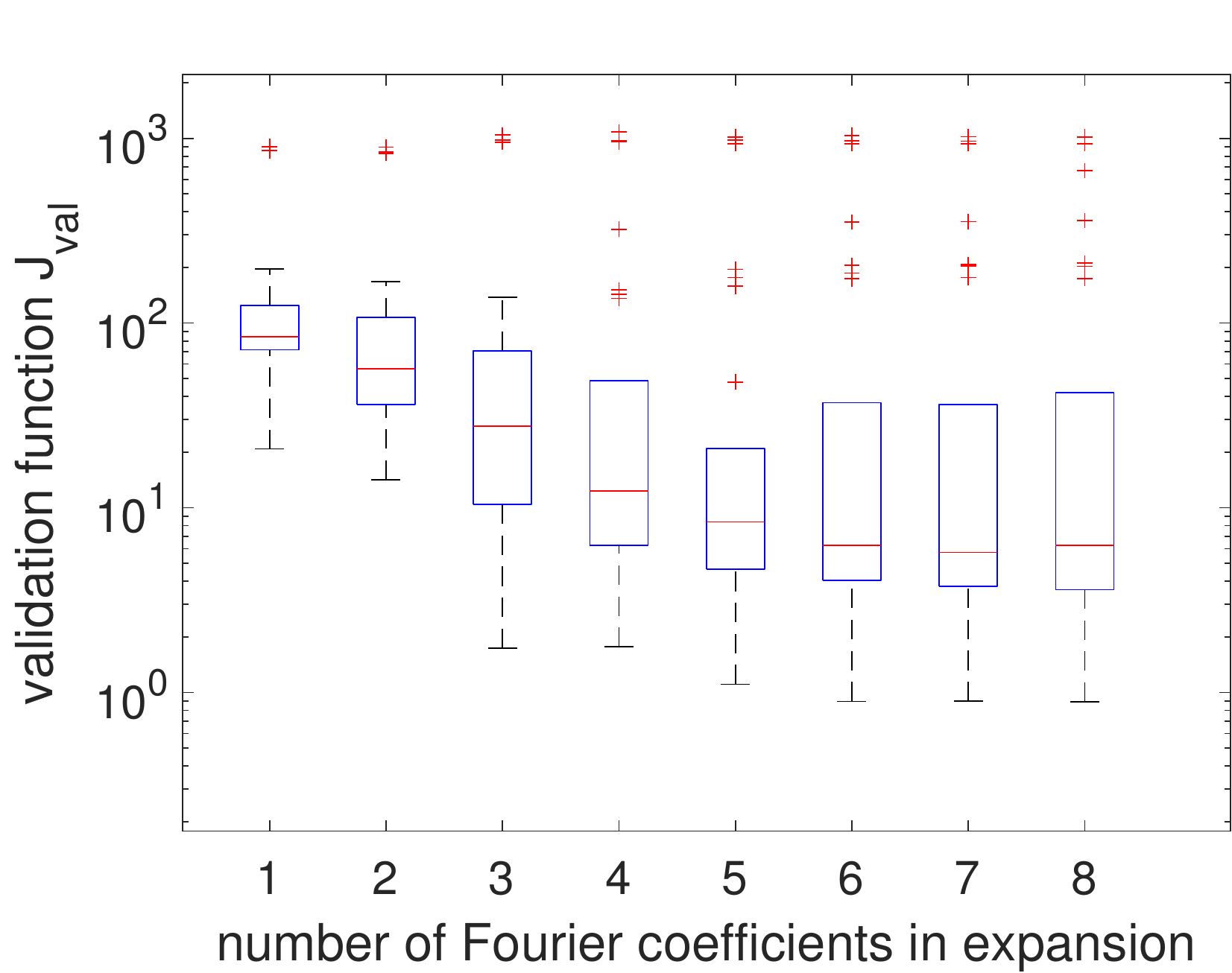}
	 	\caption{\label{FIG11}{Validation function $J_{\mathrm{val}}$ depending on number of Fourier coefficients in area function.}}
	 \end{figure}
	In the parameter fittings in the previous sections, the number of geometrical parameters $c_{i}$ and $s_{i}$ was always chosen as $M=4$. In this section, observations on varying number of parameters $M$ are discussed.\\
	Figure \ref{FIG11} illustrates the validation function $J_{\mathrm{val}}$ from \eqref{eq:Jval} computed for varying parameter $M$. The underlying data for the boxplot are the refined parameter fittings of all ten subjects and all three different impedance models (i)-(iii). We observe only small changes for $M\ge 4$. Thus, due to increasing computational cost in the following experiments $M=4$ was chosen as a reasonable compromise between accuracy and computational cost.

	\section{Validation\label{sec:5}}
	\subsection{Transfer impedance}
	\begin{figure}[H]
		\centering
		\includegraphics[width=0.55\textwidth]{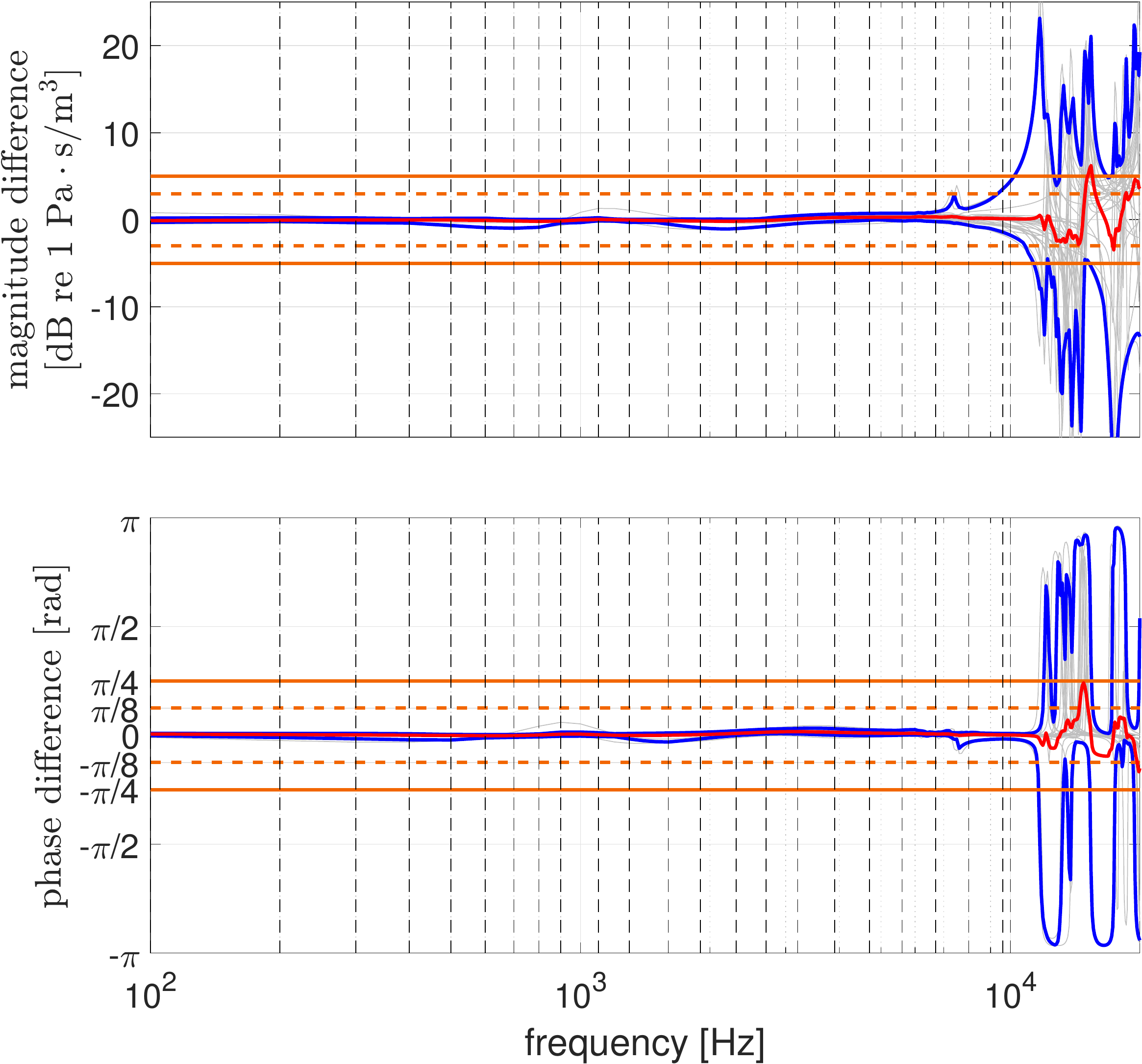}
		\caption{\label{FIG12}{Validation for the parameter fitting considering differences of the transfer impedance for frequencies for up to 10\,kHz. Straight solid line: 5\,dB or 45° difference. Straight dashed line: 3\,dB or 22.5° difference. Blue lines: 5 and 95 quantile. Red line in between: Mean. Dashed gray lines: Frequencies used during the parameter fitting.}}
	\end{figure}
	In this section validation results are presented. To this end, ten different ear canals with three different impedance models at the eardrum were considered as discussed in section \ref{sec:3}. 
		\begin{figure}[H]
		\centering
		\includegraphics[width=0.55\textwidth]{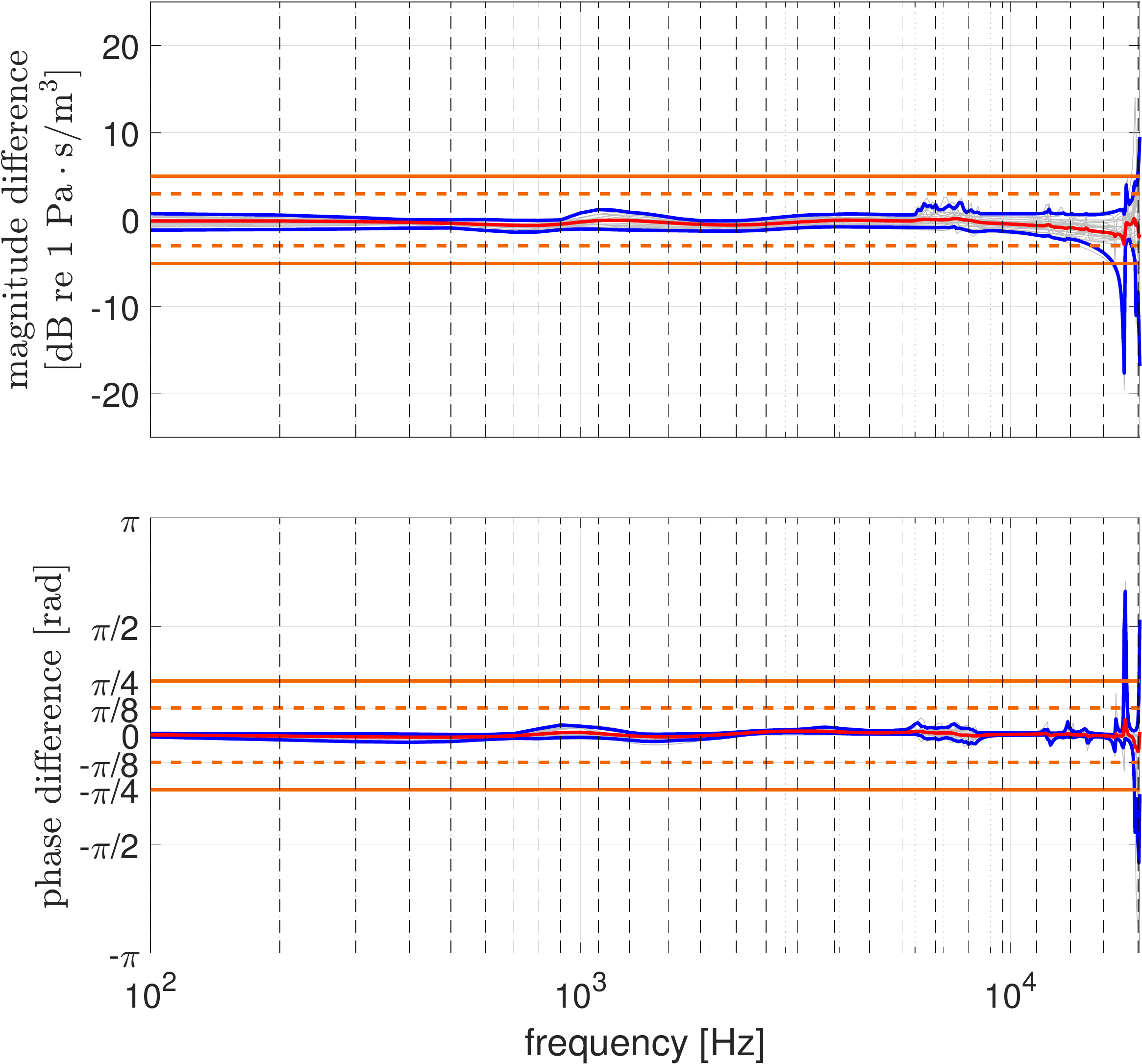}
		\caption{\label{FIG13}{Validation for the parameter fitting considering differences of the transfer impedance for frequencies for up to 20\,kHz exluding subject 9 that showed unusual behavior for high frequencies. Straight solid line: 5\,dB or 45° difference. Straight dashed line: 3\,dB or 22.5° difference. Blue lines: 5 and 95 quantile. Red line in between: Mean. Dashed gray lines: Frequencies used during the parameter fitting.}}
	\end{figure}
	In figure \ref{FIG12} the transfer impedance level and phase difference between the data and the  parameter fitting, i.e., 
	\[
	20\log_{10}\left(\left\lvert\frac{Z_{\mathrm{tr}}}{Z_{\mathrm{tr,data}}}\right\rvert\right) \text{ and } \arg\left(\frac{Z_{\mathrm{tr}}}{Z_{\mathrm{tr,data}}}\right),
	\] 
	are plotted. The blue lines indicate the 5\% and 95\% quantile over the 30 parameter fittings based on 30 data sets, i.e., 10 different subjects and three different impedance models at the eardrum, in total.

 	The parameter fitting was done for twelve initial parameter sets in all cases. Only the one with the lowest cost function over all 200 frequencies was kept. The results are illustrated in figure \ref{FIG12}. In this case the results show magnitude differences of less than 1\,dB for all frequencies in the range of up to 7\,kHz, and still magnitude differences smaller than 5\,dB for frequencies between 7\,kHz and 10\,kHz. In most cases the magnitude difference is even smaller than 3\,dB for frequencies up to 10\,kHz. For larger frequencies the differences are very large since only frequencies up to 10\,kHz were used in the parameter fitting.  Note however, that for most geometries reasonable results were also achieved for higher frequencies, see figure \ref{FIG13}.

	\subsection{Ear canal area function}
	During the  parameter fitting, the estimated parameters for the ear canal area function can be used to construct an area function to get an approximate model of the geometry. 
	\begin{figure}[H]
		\includegraphics[width=1\textwidth]{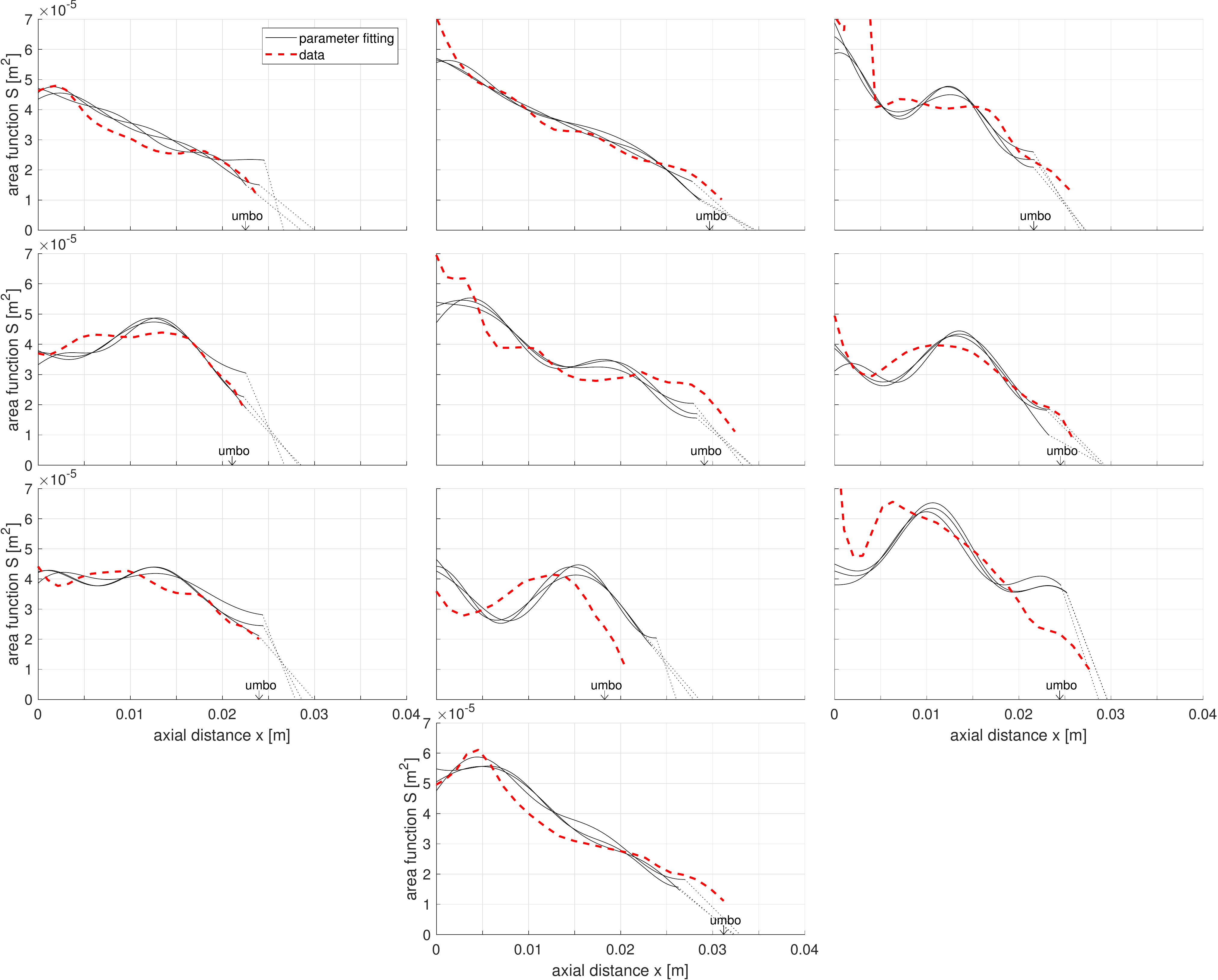}
		\caption{\label{FIG14}{Results for the ear canal area functions for all ten subjects and all three impedance models 10\,kHz.}}
	\end{figure}
	In figure \ref{FIG14} area functions taken from the parameter fitting (using logarthmically distributed frequencies up to 10\,kHz) are displayed comparing it with the area function computed from the 3D geometry using the method described in \cite{Stinson.1989}. Main features for most ear canals could be reproduced reasonably well. Note that the geometry used for the generation of the data was cut at the first bend of the center axis obtained from using the VMTK toolbox. This center axis differs from the one proposed by \cite{Stinson.1989}. Especially the entrance surface will not be perpendicular to the center axis by \cite{Stinson.1989} which explains the large deviations at the entrance of the area function for some subjects like subject 3 and subject 9. Overall, the approximations from the parameter fitting are quite close to the area function generated with the method from \cite{Stinson.1989}. The area functions are also largely unaffected by the impedance model used for the data generation.

	\subsection{Impedance at eardrum}
	\label{subsec:validationImpedance}	
	\begin{figure}[hbt!]
		\includegraphics[width=1\textwidth]{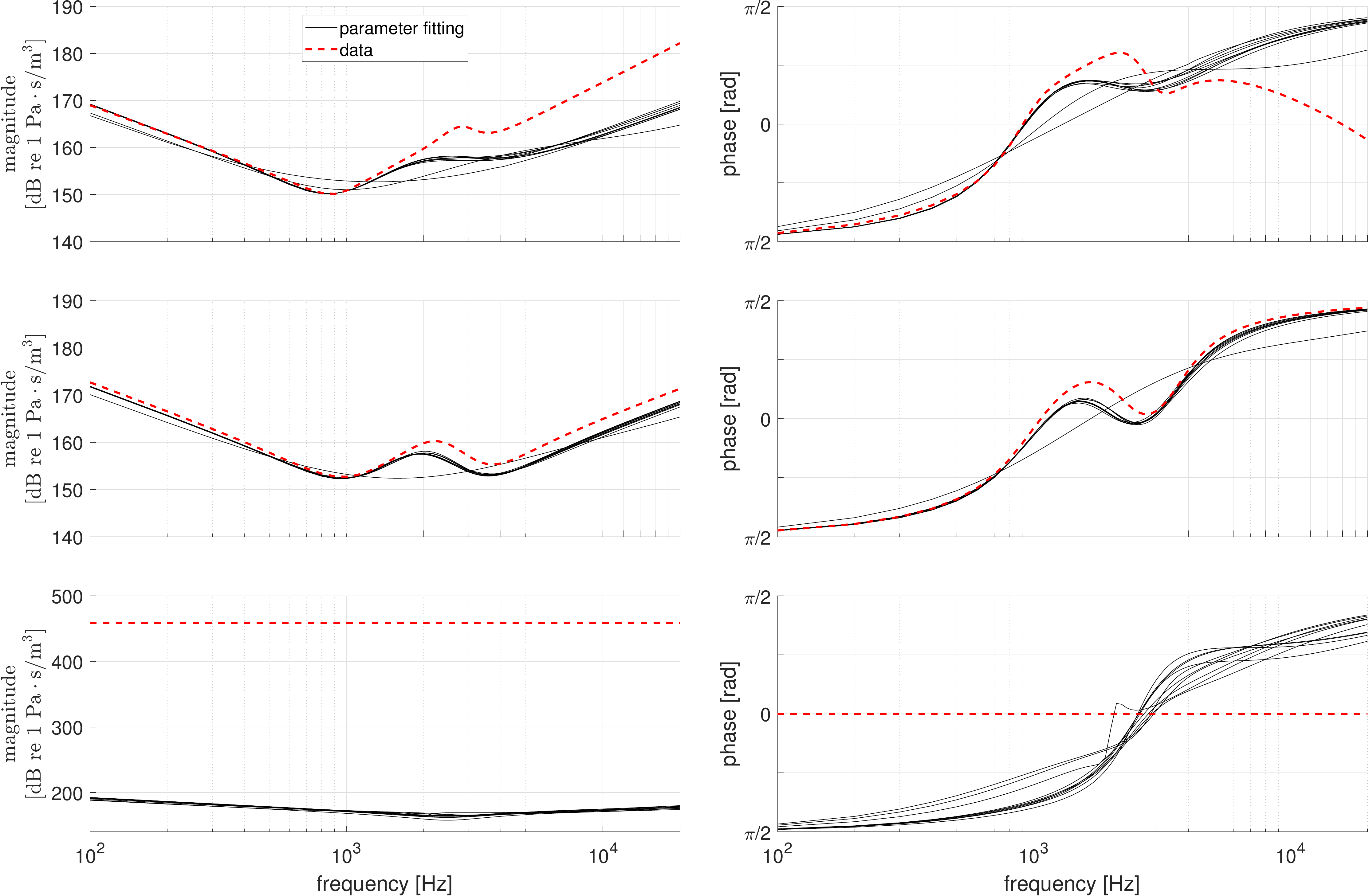}
		\caption{\label{FIG15}{Estimated impedances at the eardrum for different underlying model data: (i) Hudde (ii) two-resonator (iii) rigid.}}
	\end{figure}
	The parameter fitting also lead to an approximation of the impedance at the eardrum. This is illustrated in figure \ref{FIG15}. The data of the two-resonator model (ii) can be reconstructed well since this is the impedance model we used in the parameter fitting. The model by \cite{Hudde.1998c} is approximated well in the low frequency range. The impedance condition (iii) could not be approximated well since the bounds on the parameters do not allow impedances of such high magnitude. However, one can observe that the magnitude is larger than in the other cases, such that the resulting  parameter fitting is closer to the rigid case. In all cases, the true eardrum impedance had a larger magnitude than the estimated one at frequencies above 2\,kHz.

	\section{Conclusion\label{sec:6}}
	In this paper a method to estimate the transfer impedance from input impedance is described and validated. Using synthetically generated data the method produces accurate results, especially in the frequency range of 100\,Hz to 7\,kHz. We can furthermore estimate the ear canal area function and the eardrum impedance, where for the area function the main features of the geometry were recognized. The eardrum impedance was estimated well in the low frequency range if the underlying impedance model was fitting. For higher frequencies, however, the eardrum impedance estimation was underestimated in general.\\
	Future work should involve the validation and enhancement of the method using in vivo measurements.

	\section*{Acknowledgments}
	This research was partially funded by the Deutsche Forschungsgemeinschaft (DFG, German
	Research Foundation) – Project ID 352015383 – SFB 1330 C1
	
	\section*{Data availability statement}
	\begin{flushleft}
		\hyphenpenalty=10000
		\exhyphenpenalty=10000
			An implementation of the complete model is available on GitHub, under the reference {https://github.com/nickwulbusch/ear-canal-parameter-fitting-1d}.
	\end{flushleft}

\bibliography{sampbib}

\begin{thebibliography}{10}
\providecommand{\url}[1]{#1}
\csname url@samestyle\endcsname
\providecommand{\newblock}{\relax}
\providecommand{\bibinfo}[2]{#2}
\providecommand{\BIBentrySTDinterwordspacing}{\spaceskip=0pt\relax}
\providecommand{\BIBentryALTinterwordstretchfactor}{4}
\providecommand{\BIBentryALTinterwordspacing}{\spaceskip=\fontdimen2\font plus
\BIBentryALTinterwordstretchfactor\fontdimen3\font minus
  \fontdimen4\font\relax}
\providecommand{\BIBforeignlanguage}[2]{{%
\expandafter\ifx\csname l@#1\endcsname\relax
\typeout{** WARNING: IEEEtran.bst: No hyphenation pattern has been}%
\typeout{** loaded for the language `#1'. Using the pattern for}%
\typeout{** the default language instead.}%
\else
\language=\csname l@#1\endcsname
\fi
#2}}
\providecommand{\BIBdecl}{\relax}
\BIBdecl

\bibitem{Blau.2010}
M.~Blau, T.~Sankowsky, P.~Roeske, H.~Mojallal, M.~Teschner, and C.~Thiele,
  ``Prediction of the sound pressure at the ear drum in occluded human cadaver
  ears,'' \emph{Acta Acustica united with Acustica}, vol.~96, no.~3, pp.
  554--566, 2010.

\bibitem{Webster.1919}
A.~G. Webster, ``Acoustical impedance and the theory of horns and of the
  phonograph,'' \emph{Proceedings of the National Academy of Sciences of the
  United States of America}, vol.~5, no.~7, pp. 275--282, 1919.

\bibitem{Stinson.1989}
M.~R. Stinson and B.~W. Lawton, ``Specification of the geometry of the human
  ear canal for the prediction of sound-pressure level distribution,''
  \emph{The Journal of the Acoustical Society of America}, vol.~85, no.~6, pp.
  2492--2503, 1989.

\bibitem{Hudde.1999}
H.~Hudde, A.~Engel, and A.~Lodwig, ``Methods for estimating the sound pressure
  at the eardrum,'' \emph{The Journal of the Acoustical Society of America},
  vol. 106, no. 4 Pt 1, pp. 1977--1992, 1999.

\bibitem{Rasetshwane.2011}
D.~M. Rasetshwane and S.~T. Neely, ``Inverse solution of ear-canal area
  function from reflectance,'' \emph{The Journal of the Acoustical Society of
  America}, vol. 130, no.~6, pp. 3873--3881, 2011.

\bibitem{SankowskyRothe.2011}
T.~Sankowsky-Rothe, M.~Blau, E.~Rasumow, H.~Mojallal, M.~Teschner, and
  C.~Thiele, ``Prediction of the sound pressure at the ear drum in occluded
  human ears,'' \emph{Acta Acustica united with Acustica}, vol.~97, no.~4, pp.
  656--668, 2011.

\bibitem{SankowskyRothe.2015}
T.~Sankowsky-Rothe, M.~Blau, S.~K{\"o}hler, and A.~Stirnemann, ``Individual
  equalization of hearing aids with integrated ear canal microphones,''
  \emph{Acta Acustica united with Acustica}, vol. 101, no.~3, pp. 552--566,
  2015.

\bibitem{Souza.2014}
N.~N. Souza, S.~Dhar, S.~T. Neely, and J.~H. Siegel, ``Comparison of nine
  methods to estimate ear-canal stimulus levels,'' \emph{The Journal of the
  Acoustical Society of America}, vol. 136, no.~4, pp. 1768--1787, 2014.

\bibitem{Rasetshwane.2012}
D.~M. Rasetshwane, S.~T. Neely, J.~B. Allen, and C.~A. Shera, ``Reflectance of
  acoustic horns and solution of the inverse problem,'' \emph{The Journal of
  the Acoustical Society of America}, vol. 131, no.~3, pp. 1863--1873, 2012.

\bibitem{Nelder.1965}
J.~A. Nelder and R.~Mead, ``A simplex method for function minimization,''
  \emph{The Computer Journal}, vol.~7, no.~4, pp. 308--313, 1965.

\bibitem{Wulbusch}
\BIBentryALTinterwordspacing
N.~Wulbusch, ``Ear canal parameter fitting 1d [code],'' 2023. [Online].
  Available:
  \url{https://github.com/nickwulbusch/ear-canal-parameter-fitting-1d}
\BIBentrySTDinterwordspacing

\bibitem{A.Stirnemann.}
A.~Stirnemann, ``Ein {M}ittelohrmodell basierend auf der
  {A}ussenohr-{T}ransferimpedanz,'' \emph{Fortschritte der Akustik - DAGA 2011,
  D{\"u}sseldorf}, 2011.

\bibitem{Hudde.1998c}
H.~Hudde and A.~Engel, ``Measuring and modeling basic properties of the human
  middle ear and ear canal. part iii: Transfer functions and model
  calculations,'' \emph{Acta Acustica}, no.~84, pp. 1091--1109, 1998.

\bibitem{Rosowski.1990}
J.~J. Rosowski, P.~J. Davis, S.~N. Merchant, K.~M. Donahue, and M.~D. Coltrera,
  ``Cadaver middle ears as models for living ears: comparisons of middle ear
  input immittance,'' \emph{The Annals of otology, rhinology, and laryngology},
  vol.~99, no. 5 Pt 1, pp. 403--412, 1990.

\bibitem{Hudde.1998a}
H.~Hudde and A.~Engel, ``Measuring and modeling basic properties of the human
  middle ear and ear canal. part i: Model structure and measuring techniques,''
  \emph{Acta Acustica}, no.~84, pp. 720--738, 1998.

\bibitem{fminsearchbnd}
\BIBentryALTinterwordspacing
J.~D'Errico. fminsearchbnd, fminsearchcon. [Online]. Available:
  \url{https://www.mathworks.com/matlabcentral/fileexchange/8277-fminsearchbnd-fminsearchcon}
\BIBentrySTDinterwordspacing

\bibitem{inproceedings}
R.~Roden and M.~Blau, ``The {IHA} database of human geometries including torso,
  head and complete outer ears for acoustic research,'' 08 2020, 49th
  International Congress and Exposition on Noise Control Engineering,
  Inter-Noise, Seoul.

\bibitem{Roden.2021}
------, ``The {IHA} database of human geometries including torso, head and
  complete outer ears for acoustic research,'' 2021.

\bibitem{vmtk}
L.~Antiga, ``Patient-specific modeling of geometry and blood flow in large
  arteries,'' 01 2002, doctoral thesis, Politecnico di Milano, dipartimento di
  bioingegneria.

\bibitem{matlabCurve}
\BIBentryALTinterwordspacing
A.~Mjaavatten, ``Curvature of a 1d curve in a 2d or 3d space.'' [Online].
  Available:
  \url{https://www.mathworks.com/matlabcentral/fileexchange/69452-curvature-of-a-1d-curve-in-a-2d-or-3d-space,
  (last visited March 8, 2023).}
\BIBentrySTDinterwordspacing

\bibitem{frenet}
\BIBentryALTinterwordspacing
D.~Claxton, ``\BIBforeignlanguage{MATLAB Central File Exchange}{Frenet}.''
  [Online]. Available:
  \url{https://www.mathworks.com/matlabcentral/fileexchange/11169-frenet, (last
  visited March 8, 2023).}
\BIBentrySTDinterwordspacing

\bibitem{comsol}
\BIBentryALTinterwordspacing
COMSOL, ``Comsol multiphysics \textregistered,'' \emph{Stockholm, Sweden}.
  [Online]. Available: \url{www.comsol.com}
\BIBentrySTDinterwordspacing

\end{thebibliography}
\end{document}